\magnification1170
\input graphicx
\input  amssym.tex
\overfullrule0pt

\def\sqr#1#2{{\vcenter{\hrule height.#2pt              
     \hbox{\vrule width.#2pt height#1pt\kern#1pt
     \vrule width.#2pt}
     \hrule height.#2pt}}}
\def\square{\mathchoice\sqr{5.5}4\sqr{5.0}4\sqr{4.8}3\sqr{4.8}3}
\def\qed{\hskip4pt plus1fill\ $\square$\par\medbreak}

\centerline{\bf Julia sets for polynomial diffeomorphisms of ${\Bbb C}^2$ are not semianalytic}


\medskip\centerline{Eric Bedford and Kyounghee Kim}

\bigskip\noindent{\bf \S0.  Introduction. }  If $X$ is a complex manifold, and $f:X\to X$ is a holomorphic mapping, then the Fatou set is the open set where the iterates $f^n:=f\circ\cdots\circ f$ are locally equicontinuous.  Equivalently, these are the points where $f$ is Lyapunov stable.  The complement of the Fatou set is the Julia set.  While we refer to this as {\it the} Julia set,  it is possible to define several Julia sets $J_k$, (see [FS1] and [U2]), in which case the complement of the Fatou set is the first Julia set $J_1$.  In dimension 1, the principal case is where $X={\Bbb P}^1$ is the Riemann sphere, and $f$ is a rational function.  In this case, Fatou showed that if $J$ has a tangent at some point, then $J$ is either a circle or a circular arc.  In the case of the circle, $f$ is conjugate to $z^d$ for $d\in{\Bbb Z}$, $|d|\ge2$; and in the case of an arc, $f$ is conjugate to a Chebyshev polynomial.  In higher dimension, there are of course product maps, and in this case the Julia set is a product.  There are also nontrivial examples of polynomial maps for which the Julia set is (real) algebraic; examples were given in ${\Bbb C}^2$ by Nakane  [N] and in ${\Bbb C}^3$ by Uchimura [Uch1--3].  

These maps discussed above are non-invertible; in the sequel we consider invertible maps.  In this case, we have both a forward Julia set $J^+:=J(f)$ and a backward Julia set $J^-:=J(f^{-1})$.  
The invertible polynomial maps of ${\Bbb C}^2$ have been classified by Friedland and Milnor [FM].  The polynomial diffeomorphisms with nontrivial dynamical behavior are compositions of generalized H\'enon maps, and each such composition has a degree $d$.  (See [BS1], [FS], [HO1] for the basic dynamical properties of these maps.)   By [FM] and [S] it follows that the topological entropy of $f$ is $\log(d)$.  Hubbard [H] defined the escape locus $U^+$ for such a map $f$, and it is easily seen that $J^+=\partial U^+$.  By [BS], $J^+$ cannot contain an algebraic curve, so it follows (see Proposition 1.2) that: {\sl Neither $J^+$ nor $J^-$ can be (real) algebraic.}  Our main result is:
\proclaim Theorem.  Let $f$ be a polynomial diffeomorphism of ${\Bbb C}^2$ with positive entropy.  Then neither $J^+$ nor $J^-$ is a semianalytic set.

Forn\ae ss and Sibony [FS2] showed that, for a generic H\'enon map, the Julia set is neither smooth nor semianalytic.  In [BK] we showed that $J^+$ can never be $C^1$ smooth.   However, the Julia sets in [N] and [Uch1--3] have singular points and thus are not $C^1$, and this non{-}smoothness was our motivation for the present Theorem.
%

\bigskip\noindent{\bf \S1.  Levi flat hypersurfaces. }  Let $U\subset{\Bbb C}^k$ be an open subset.  A function $\rho$ on $U$ is said to be {\it real analytic} if for every $q\in U$, $\rho$ can be written as a real power series which converges in a neighborhood of $q$.  Let us suppose that $q=0$ and write
$$\rho(z,\bar z)=\sum_{I,J} c_{I,J}z^I\bar z^J$$
where $I=(i_1,\dots,i_k)$ is a $k$-tuple of  nonnegative integers, and $z^I=z_1^{i_1}\cdots z_k^{i_k}$, and similarly for $J$ and $\bar z^J$.  We may treat $z$ and $\bar z$ as independent variables and write
$$\rho(z,\bar w) = \sum_{I,J} c_{I,J}z^I\bar w^J$$
The reality condition on $\rho$ is that $c_{I,J}=\overline{ c_{J,I}}$, which means that $\rho(z,\bar w)=\overline{\rho(w,\bar z)}$.   A set $X$ is {\it real analytic} if it can be written locally as $X\cap U=\{\rho=0\}$.  A point $x_0\in X$ is said to be {\it regular} if $X$ is a smooth manifold in a neighborhood of $x_0$.  We write $Reg(X)$ for the set of regular points, and $Reg(X)$ is dense in $X$ (see [BM]), although the dimension may be different at different points.

A real hypersurface is said to be {\it Levi flat} if it is (locally) pseudoconvex from both sides.  Recall that $G^+$ is pluriharmonic on the set $\{G^+>0\}$, so: { \sl If the set $J^+=\partial\{G^+>0\}$ is $C^1$ smooth at some point, then it is Levi flat there. }  A real analytic set is said to be Levi flat if it is Levi flat at each regular point.   If $X$ is real analytic, Levi flat hypersurface, then at each regular point, there is a local holomorphic coordinate system such that $X$ is locally given as $\{z_1 + \bar z_1=0\}$.  At singular points, the situation is more complicated.


  A semianalytic set is given locally as a finite, disjoint union of sets of the form $\{r_j=0, s_j>0\}$, where $r_j$ and $s_j$ are real analytic.  (See Bierstone and Milman [BM] for further information on semianalyticity.)  The following Lemma allows us to replace a semianalytic set $J^+$ by an analytic set $X$.

\proclaim Lemma 1.1.  Suppose that $J^+$ is semianalytic, and $p\in J^+$ is fixed under $f$.  Then there is a neighborhood $U$ of $p$ and a real analytic Levi flat hypersurface $X:=\{\rho=0\}\subset U$ containing $p$ such that:  
\item{(i)} $X\cap V=J^+\cap V\ne\emptyset$ for some open set $V\subset U$,
\item{(ii)} $X$ is locally invariant in the sense that $f(X)\cap U\subset X$.

\noindent{\it Proof. }  If $J^+$ is semianalytic in a neighborhood of $p$, then there is an open set $U$ containing $p$ such that $J^+\cap U$ is a finite union of sets $X_j=\{r_j=0, s_j>0\}$, where $r_j, s_j$ are real analytic on $U$.   We may assume that $p$ belongs to the closure of each $X_j$.  Passing to an iterate of $f$, we may assume that each $X_j$ is invariant, in the sense that  $f(X_j\cap U)\subset X_j$.  Since the set of regular points is dense, we may suppose that at least one of the $X_j$ intersects $Reg(J^+)$.  Let us set $\rho:=r_j$ and $X:=\{\rho=0\}$.  We may assume that $X$ is an irreducible subvariety of $U$; otherwise, pass to a component.  We have local invariance (ii).  Since $X$ intersects $Reg(J^+)$, we may choose an open set $V$ such that condition (i) holds, and such that $V\cap J^+$ consists of regular points.  Since $J^+\cap V$ is regular, it is Levi flat, and since $X$ is irreducible, it is Levi flat, too.  \qed


Let us discuss the hypersurface $X=\{\rho=0\}$, where $\rho(z,\bar w)$ converges for $z,w\in U$.   If for fixed $w\in U$, $\rho(z,\bar w)= 0$ for all $z$, we say that $X$ is {\it Segr\`e degenerate} at $w$.  If  $X$ is not degenerate at $w\in U$, then the Segr\`e variety, which is defined as
$$Q_w:=\{z\in U: \rho(z,\bar w)=0\},$$
is a proper subvariety of $U$.  (In other words, the condition that $w$ is Segr\`e degenerate means that the Segr\`e variety is the whole open set $U$.)  We may choose the defining function $\rho$ to be minimal, which means that if $\rho'$ is any other defining function, then $\rho'=h\rho$.  The family of Segr\`e varieties is independent of the choice of minimal defining function.  

At this stage, we can conclude that $J^\pm $ cannot be algebraic.  
\proclaim Proposition 1.2.  Let $f$ be a polynomial diffeomorphism of ${\Bbb C}^2$ with positive entropy.  Then neither $J^+$ nor $J^-$ is an algebraic set.

\noindent{\it Proof. }  Let us suppose that $J^+=\{\rho(z,\bar z)=0\}$ is defined by a real polynomial.  A regular point $w\in J^+$ is Segr\`e nondegenerate, so  
$Q_w$ is a proper subvariety of ${\Bbb C}^2$, which is contained in $J^+$.  On the other hand, this is not possible, since by [BS] there is no subvariety of ${\Bbb C}^2$ which is contained in $K^+$.
\qed

The set  of Segr\`e degenerate points is a complex subvariety of codimension at least 2 (see [PSS, \S3]).  Thus in ${\Bbb C}^2$, the Segr\`e degenerate points are isolated, so we may assume that $U$ is sufficiently small that all points of $X\cap U-\{p\}$ are Segr\`e nondegenerate.

A basic result (see Pinchuk, Shafikov and Sukhov [PSS]) is that if $X$ is Levi flat, then for $w\in X$, then each proper Segr\`e variety $Q_w$ is contained in $X$.   We say that $p$ is {\it dicritical } if there are infinitely many distinct varieties $Q_q$ passing through $p$.  Since $X$ is irreducible, it follows that if infinitely many varieties $Q_q$ contain $p$, then all varieties $Q_q$ contain $p$.  We will make use of the following result:

\proclaim Theorem 1.3 [PSS].  A point is Segr\`e degenerate if and only if it is dicritical.

\proclaim Lemma 1.4.    Suppose that $J^+$ is semianalytic, and $p\in J^+$ is fixed.  If $X$ is as in Lemma~1.1, then $p$ is not dicritical.

\noindent{\it Proof.}   If $r_0$ is a saddle point, then by [BS2], the stable manifold $W^s(r_0)$ is dense in $J^+$.  Since there are infinitely many saddle points, we may suppose that $r_0\ne p$.   Let $q\in W^s(r_0)\cap X-\{p\}$ be a regular point of $X$.  We may assume that $q$ is Segr\`e nondegenerate, so that $Q_q$ is a complex subvariety of $X$.  Further, since the leaves of the complex foliation of a Levi flat hypersurface are unique, it follows that $W^s(r_0)$ intersects $Q_q$ in an open set.   If $p$ is dicritical, then $p\in Q_q$.  On the other hand, since $p$ is fixed, it cannot belong to $W^s(r_0)$.   Thus $\hat W^s(r_0):=W^s(r_0)\cup Q_q$ is a complex manifold which is strictly larger than $W^s(q_0)$.  (Note that we may desingularize $\hat W^s(r_0)$ if $p$ is a singular point of $Q_q$.)  Now recall that $W^s(r_0)$ is uniformized by ${\Bbb C}$, and the only Riemann surface which strictly contains ${\Bbb C}$ is the Riemann sphere, which is compact.  Since ${\Bbb C}^2$ can contain no compact, complex manifolds, we have a contradiction, by which we conclude that $Q_q$ cannot contain $p$.  Thus $p$ is not dicritical.  
\qed

Let us suppose that the multipliers of $Df$ at $p$ are $|\alpha|<|\beta|$, with $|\alpha|<1$.  Then the {\it strong stable set} of $p$ corresponding to the multiplier $\alpha$ is defined as
$$W^{ss}(p)=\{p\}\cup \{q\in {\Bbb C}^2:\lim_{n\to\infty}{1\over n}\log({\rm dist}(f^n(p),f^n(q)))=\log|\alpha|\}$$
By the Strong Stable Manifold Theorem, $W^{ss}(p)$ is a complex submanifold of ${\Bbb C}^2$ which is uniformized by ${\Bbb C}$.
The {\it local strong stable manifold} is defined as 
$$W^{ss}_{\rm loc}(p):=\{q\in W^{ss}(p): f^n(q)\in U {\rm\ for\ all\ }n\ge0\}.$$    
Let us choose coordinates $(x,y)$ near $p=(0,0)$ so that the coordinate axes are the eigenspaces for $Df(p)$.  Then if we take $U=\{|x|<r_1,|y|<r_2\}$ to be a small bidisk, then $W^{ss}_{\rm loc}(p)$ is the connected component of $W^{ss}(p)\cap U$ which contains $p$.

We conclude this section by showing that the local strong stable manifold coincides with the Segr\`e variety through $p$.  
\proclaim Lemma 1.5.  Suppose that $J^+$ is semianalytic, and $p\in J^+$ is fixed.   If $X$ is as in Lemma~1.1, then  $Q_p=W^{ss}_{\rm loc}(p)$, and the multipliers of $Df$ at $p$ are $|\alpha|<1\le|\beta|$.

\noindent{\it Proof. }  
The Jacobian determinant of $f$ is a nonzero constant $\delta$.  As was shown in [BK], if $J^+$ is semianalytic, we must have $|\delta|\le1$.  Let $\alpha$, $\beta$ be the multipliers of $Df$ at $p$.  Since $|\alpha\beta|=|\delta|<1$, we may suppose that $|\alpha|<1$.  If $|\beta|<1$, then $p$ is an attracting fixed point, which means that $p$ belongs to the interior of $K^+$.  Since $p\in J^+=\partial K^+$, we must have $|\beta|\ge1$.  Thus the eigenvalues are distinct, and we may diagonalize $Df(p)$.   We may suppose that $p=(0,0)$, and $f(x,y)=(x_1,y_1)=(\beta x + \cdots,\alpha y+\cdots)$.  Further, we may choose local coordinates such that $W^{ss}_{\rm loc}(p)=\{x=0\}$.  

If $V$ be an irreducible component of $Q_p$, and $V$ is not the same as $W^{ss}_{\rm loc}(p)$, then we may choose $U$ sufficiently small that $Q_p\cap W^{ss}_{\rm loc}(p)=Q_p\cap\{x=0\}=\{(0,0)\}=\{p\}$.  Thus for some positive integer $\mu$ we may choose a root $x^{1/\mu}$ and represent $V$ locally as a Puiseux expansion $V=\{y=\sum_{j=1}^\infty a_j x^{j/\mu} \}$.  The local invariance of $V$ at $p=(0,0)$ means that we will have $y_1=\sum_{j=1}^\infty a_j x_1^{j/\mu}$.  If $a_{j_0}$ is the first nonvanishing coefficient, we must have $\alpha = \beta^{j_0/\mu}$.  But this is impossible since $j_0/\mu>0$, and  $|\alpha|<1\le|\beta|$.  It follows, then that the only irreducible component of $Q_p$ is $\{x=0\}$.   \qed

\bigskip\noindent{\bf \S2.  Multipliers at a fixed point. }

\proclaim Lemma 2.1.  Suppose that $J^+$ is semianalytic, and $p\in J^+$ is fixed.  Then $p$ is a saddle point.

\noindent{\it Proof.}  By Proposition 1.5, we know that the multipliers of $Df$ at $p$ are $|\alpha|<1$ and $|\beta|\ge1$.  We must show that $|\beta|>1$.  If not, then $|\beta|=1$.  First, we observe that $\beta$ cannot be a root of unity.  For in that case, $p$ is a semi-attracting, semi-parabolic fixed point.  By Ueda [U1] and Hakim [H], there is a semi-parabolic basin ${\cal B}$ with $\partial{\cal B}=J^+$.  However, the boundary of ${\cal B}$ has a fractal ``cusp'' at $p$ (reminiscent of the cauliflower Julia set) and is not semianalytic.  We conclude that $\beta^k\ne1$ for all nonzero integers $k$.

Now let us use coordinates from the proof of Lemma 1.5.   Since $Q_{(0,0)}=\{x=0\}$, we may write $\rho(x,y,0,0)=x^ku(x,y)$, where $u(x,y)$ is a holomorphic function with $u(0,0)=1$.  This means that 
$$\rho(x,y,\bar x, \bar y) = x^ku(x,y) + \bar x^k \overline{u(x,y)} + \Psi(x,y,\bar x, \bar y)$$
where in the expansion  of $\rho$, all of the purely holomorphic terms are contained in $x^ku(x,y)$, and $x^k$ is the purely holomorphic part of lowest order.  Now there is a real analytic unit $h(x,y,\bar x,\bar y)$ such that $\rho\circ f = h\,\rho$, and $h(0,0)=c\ne0$ is real.  Thus the purely holomorphic part of lowest order are $cx^k$.  On the other hand, as in the proof of Lemma 1.5, we have
$$\eqalign{\rho(f(x,& y))=  \rho(x_1,y_1,\bar x_1,\bar y_1) = \cr
&=\rho(\beta x + \cdots,\alpha y+\cdots,\bar\beta\bar x + \cdots,\bar\alpha\bar y + \cdots) = \beta^k x^k + \bar\beta^k\bar x^k + \Psi_1}$$
Thus we see that the purely holomorphic terms of lowest order are $\beta^k x^k$, from which we conclude that $\beta^k$ is real, which is a contradiction. \qed

If $p$ is a saddle point, we let $W^u(p)$ be the unstable manifold at $p$.  If the multipliers are $|\alpha|<1<|\beta|$, then there is a holomorphic uniformization $\psi_p:{\Bbb C}\to W^u(p)\subset{\Bbb C}^2$ such that $\psi_p(0)=0$, and $\psi_p(\beta\zeta)=f(\psi_p(\zeta))$.  We set $J_p:=\psi_p^{-1}(J^+\cap W^u(p))$.    By the invariance of $J^+$ it follows that $J_p$ is invariant under $\zeta\mapsto\beta\zeta$.   

\proclaim Lemma 2.2.  Suppose that $J^+$ is semianalytic, and $p\in J^+$ is fixed.  Then $\beta\in{\Bbb R}$, and  $J_p$ is a straight line in ${\Bbb C}$ passing through the origin.

\noindent{\it Proof.}  $J_p$ is a semianalytic set of dimension one in ${\Bbb C}$.  Thus it has a tangent cone at the origin.  Since it is invariant under $\zeta\mapsto\beta\zeta$, we conclude that $\beta\in{\Bbb R}$, and $J_p$ consists of a finite number of (infinite) rays emanating from the origin.  We must show that there are exactly two rays whose union forms a line passing through the origin.  We consider the points $r\in J_p$ which correspond to transverse intersections between $W^s(p)$ and $W^u(p)$.  By [BLS4] this set is dense in $J_p$.   Let $\Delta_0\subset{\Bbb C}$ denote a small disk about the origin, and let $\Delta\subset{\Bbb C}$ denote a disk about $r$, small enough that $J_p\cap \Delta$ is a single segment $I$ which divides $\Delta$ into halves $\Delta'$ and $\Delta''$.

Consider the complex disks in ${\Bbb C}^2$ given by ${\cal D}_0:=\psi_p(\Delta_0)$ and ${\cal D}:=\psi_p(\Delta)$.  Since ${\cal D}$ is transverse to $W^s(p)$ at $\psi_p(r)$, we may apply the Lambda Lemma to conclude that there are disks ${\cal D}_j\subset f^j({\cal D})$ containing $f^j(\psi_p(r))$ which may be written as graphs over ${\cal D}_0$, and ${\cal D}_j\to {\cal D}_0$ in the $C^1$ topology.   Let $\gamma_j:=f^j(\psi_p(I))\cap {\cal D}_j$.  This is a smooth curve which divides ${\cal D}_j$ into halves ${\cal D}_j'$ and ${\cal D}_j''$, corresponding to the partition $\Delta=\Delta'\cup I\cup \Delta''$.  It follows that the $\gamma_j$ converge to a smooth curve $\gamma_0\subset{\cal D}_0$.  Now the Green function $G^+$ cannot vanish on all of ${\cal D}_j$, so it must be strictly positive everywhere on, say, ${\cal D}_j'$.  Since $G^+$ is continuous, we have $G^+|_{{\cal D}_j'}\to G^+|_{{\cal D}_0'}$.  If this limit of positive harmonic functions vanishes at a point of ${\cal D}_0'$, then it must vanish everywhere on  ${\cal D}_0'$.  We know that $G^+$ cannot vanish everywhere on ${\cal D}_0$, so we conclude that it must be strictly positive everythere on either ${\cal D}_0'$ or ${\cal D}_0''$.  We conclude that $\gamma_0\subset J^+\cap {\cal D}_0$.   Thus we conclude that $J_p$ contains a line, and any additional rays of $J_p$ must lie on the other side, $\Delta_0''$.  However, if there is a ray inside $\Delta''_0$, there are points of $\Delta''_0$ where $G^+>0$, and we may apply the same argument to ${\cal D}''_j$ to conclude that all of $J_p$ must lie inside the line.   \qed

\proclaim Lemma 2.3.  There is a dense set of complex lines $L\subset {\Bbb C}^2$ such that $K^+\cap L$ contains interior.

\noindent {\it Proof.}   If $L\subset{\Bbb C}^2$ is a complex line, then by [FM], $L\cap J^+$ is compact.  Since $J^+$ is semianalytic of pure dimension 3, it follows that for generic $L$, $J^+\cap L$ has real dimension $\le 1$.  In fact, it has pure dimension 1 since there can be no component of dimension zero (which would be an isolated point) because  $J^+$ is the boundary of $\{G^+>0\}$.  Now let us fix a complex line $L_0$ such that  $J^+\cap L_0$ consists of a finite number of segments and closed curves.  We will show that there exists a line $L$ arbitrarily close to $L_0$  such that $K^+\cap L$ contains interior.

Otherwise, for all $L$ in a neighborhood of $L_0$, $J^+\cap L$ is simply connected, and thus it must be a tree.   We let ${\cal E}_L\ne\emptyset$ denote the set of endpoints of this tree.  Now recall that in a neighborhood of $L_0\cap J^+$, $J^+$ is a finite union of disjoint strata $M_j=\{r_j=0,s_j>0\}$.  If $M_j$ has dimension 3, then for generic $L_0$, $M_j\cap L$ will have dimension 1 for all $L$ near $L_0$.   Thus the endpoints ${\cal E}_L$ can come only from $M_j$ of dimension 2.   Since $J^+$ is Levi flat,  $M_j$ must be complex analytic.  Further, since $J^+$ is Levi flat, we may follow it globally to conclude that $M$ is contained inside $\tilde M$, which is a local variety inside the boundary of $J^+$.   Thus $\tilde M\subset J^+$ is a subvariety of ${\Bbb C}^2$.  However, there is no complex subvariety contained in $K^+$ (see [BS] or [FS1]), which is a contradiction. 

We conclude from the contradiction that  for some $L$ close to $L_0$, $J^+\cap L$ is not simply connected and thus divides $L$ into (at least) two connected components.  Only one of these components can be unbounded, so we let $\omega\subset L$ denote a bounded component of the complement of $J^+\cap L$.  On the other hand, $G^+\ge0$ vanishes on $J^+$, so by the maximum principle, $G^+=0$ on $\omega$, so $\omega\subset K^+=\{G^+=0\}$. \qed

%
%
%

\proclaim Lemma 2.4.   Let $f$ be a polynomial diffeomorphism of ${\Bbb C}^2$ with positive entropy, and let $d$ be the degree of $f$.   If $J^+$ is a semianalytic set, and if $p\in J^+$ is a fixed point, then $d$ is one of the eigenvalues of $Df$ at $p$.

\noindent {\it Proof }  We continue to let $\psi_p:{\Bbb C}\to W^u(p)$ denote the uniformization of $W^u(p)$, and we define $g(\zeta):=G^+(\psi_p(\zeta))$, which is subharmonic on ${\Bbb C}$ and satisfies the functional equation $g(\beta\zeta)=d\cdot g(\zeta)$.  By Lemma 2.2, we may assume that $J_p$ is the real axis. Thus on the upper/lower half plane, $g(\zeta)=c^\pm \Im (\zeta)$ for some constants $c^+\ge0$ and $c^-\le0$, which are not both zero.  By the functional equation, we have $c^+\Im (\beta\zeta)=d\cdot c^+\Im (\zeta)$ if $\beta>0$, so $\beta=d$ in this case.  If $\beta<0$, then we have $c^+=-c^-$, and $\beta=-d$.  

Now we will show that one of the $c^\pm$ is zero, so we must have $\beta=d$.   By Lemma 2.3, we may choose a  $L\subset{\Bbb C}^2$ such that $K^+\cap L$ contains an interior component $\omega$.  We may choose a point $r\in W^s(p)\cap \partial\omega$ which is a regular point of $\partial\omega$.  Further, we may suppose that $L$ is transverse to $W^s(p)$ at $r$.  Now we let $\Delta\subset L$ denote a small disk containing $r$, so that $\Delta\cap \partial\omega$ is a smooth arc which divides $\Delta$ into two open components.  We have $G^+=0$ on $\omega\cap \Delta$ and $G^+>0$ on the complementary component.   Now we apply the Lambda Lemma as we did in Lemma 2.2, and we conclude that $G^+=0$ on one of the components of the complement of ${\cal D}_0\cap J^+\subset W^u(p)$.  Thus we have $c^+=0$ or $c^-=0$.  \qed

Our Theorem is a consequence of Lemma 2.4:
\medskip\noindent{\it Proof of the Theorem. }   If $J^+$ is semianalytic, then by [BK], $f$ must be dissipative (volume decreasing).  Then by [BS2], there can be at most one fixed point $p\in {\rm int}(K^+)$.  Thus every fixed point, with at most one exception, is contained in $J^+$.  By Lemma 2.4, $d$ is a multiplier for $Df$ at each fixed point, except possibly one.  However,  by Proposition 5.1 of [BK], this is not possible, so $J^+$ cannot be semianalytic.      \qed

\bigskip\centerline{\bf References}
\medskip
\item{[BK]}  E. Bedford and K. Kim,  No smooth Julia sets for polynomial diffeomorphisms of ${\Bbb C}^2$ with positive entropy, J. of Geometric Analysis, to appear.

\item{[BLS4]} E. Bedford, M. Lyubich and J. Smillie,  Polynomial diffeomorphisms of ${\bf C}^2$. IV. The measure of maximal entropy and laminar currents. Invent. Math. 112 (1993), no. 1, 77--125.

\item{[BS]}  E. Bedford and J. Smillie, Fatou-Bieberbach domains arising from polynomial automorphisms. Indiana Univ. Math. J. 40 (1991), no. 2, 789--792. 

\item{[BS1]}  E. Bedford and J. Smillie, Polynomial diffeomorphisms of ${\bf C}^2$: currents, equilibrium measure and hyperbolicity. Invent. Math. 103 (1991), no. 1, 69--99. 

\item{[BS2]} E. Bedford and J. Smillie, Polynomial diffeomorphisms of ${\bf C}^2$. II. Stable manifolds and recurrence. J. Amer. Math. Soc. 4 (1991), no. 4, 657--679. 

\item{[BM]}  E. Bierstone and P. Milman,  Semianalytic and subanalytic sets,   Inst.\ Hautes \'Etudes Sci. Publ. Math. No. 67 (1988), p. 5--42.

\item{[FS1]} J.E. Forn\ae ss and N. Sibony, Complex dynamics in higher dimensions. Notes partially written by Estela A. Gavosto. NATO Adv. Sci. Inst. Ser. C Math. Phys. Sci., 439, Complex potential theory (Montreal, PQ, 1993), 131--186, Kluwer Acad. Publ., Dordrecht, 1994.

\item{[FS2]} J.E. Forn\ae ss and N. Sibony,  Complex H\'enon mappings in ${\bf C}^2$ and Fatou-Bieberbach domains. Duke Math. J. 65 (1992), no. 2, 345--380. 

\item{[FM]}  S. Friedland and J. Milnor,  Dynamical properties of plane polynomial automorphisms.
Ergodic Theory Dynam. Systems 9 (1989), no. 1, 67--99. 

\item{[Ha]}  M. Hakim,  Attracting domains for semi-attractive transformations of ${\Bbb C}^p$. Publ. Mat. 38 (1994), no. 2, 479--499.

\item{[H]}  J.H. Hubbard,  The H\'enon mapping in the complex domain. Chaotic dynamics and fractals (Atlanta, Ga., 1985), 101--111, Notes Rep. Math. Sci. Engrg., 2, Academic Press, Orlando, FL, 1986.

\item{[HO1]}  J.H. Hubbard and R. Oberste-Vorth,  H\'enon mappings in the complex domain. I. The global topology of dynamical space. Inst.\ Hautes \'Etudes Sci. Publ. Math. No. 79 (1994), 5--46.

\item{[N]}  S. Nakane,   External rays for polynomial maps of two variables associated with Chebyshev maps. J. Math. Anal. Appl. 338 (2008), no. 1, 552--562.

\item{[PSS]}  S. Pinchuk, R. Shafikov, A. Sukhov,  Dicritical singularities and laminar currents on Levi-flat hypersurfaces,   arXiv:1606.02140



\item{[S]}  J. Smillie,  The entropy of polynomial diffeomorphisms of ${\Bbb C}^2$. Ergodic Theory Dynam. Systems 10 (1990), no. 4, 823--827. 

\item{[Uch1]}  K. Uchimura, The dynamical systems associated with Chebyshev polynomials in two variables, Internat. J. Bifur. Chaos Appl. Sci. Engrg. 6, 12b(1996) 2611--2618.

\item{[Uch2]}  K. Uchimura,  The sets of points with bounded orbits for generalized Chebyshev mappings, Internat. J. Bifur. Chaos, Vol. 1, No. 1 (2001) 91--107. 

\item{[Uch3]}  K. Uchimura,  Holomorphic endomorphisms of ${\Bbb P}^3({\Bbb C})$ related to a Lie algebra of type $A_3$ and catastrophe theory, Kyoto J. Math. vol 57 no 1 (2017) 197--232.

\item{[U]} T. Ueda,  Local structure of analytic transformations of two complex variables.~I.
{\sl J. Math. Kyoto Univ.} 26 (1986), no. 2, 233--261. 

\item{[U2]}  T. Ueda, Fatou sets in complex dynamics on projective spaces. J. Math.\ Soc.\ Japan 46 (1994), no. 3, 545--555.

\bigskip
\rightline{Eric Bedford}

\rightline{Stony Brook University}

\rightline{Stony Brook, NY 11794}

\rightline{\tt ebedford@math.stonybrook.edu}

\bigskip
\rightline{Kyounghee Kim}

\rightline{Florida State University}

\rightline{Tallahassee, FL 32306}

\rightline{\tt kim@math.fsu.edu}

\bye